# Asymptotic expansion of a function defined by power series

Mihail M. Nikitin

Sevastopol National Technical University,
Sevastopol, Ukraine, 99053

**Abstract**

We present a sufficient condition of existence of asymptotic expansion in negative power series for a function defined by Taylor series and unitary formulas for coefficients of this expansion. An example of computing scheme for arctangent function is represented.

**Key words and phrases**

Power series analysis, asymptotic expansion in prescribed system of functions, Padé approximations, quadratic approximations.

## Introduction

In accordance with the theory of analytic continuation of complex variable function all properties of a function analytic in some point are defined by its power expansion in this point [1]. The central practical problem of the theory is research of function properties immediately on series coefficients prescribed.

At present there exist several methods of this problem solution based on assumption that we possess beforehand some information either about class of the function or about general tend of its power expansion terms behaviour [2].

There are: method of Padé approximations [2], method of $G^3 J$ - approximations [3,4], method of quadratic approximations [5], method of Levin approximations [6]. Also methods similar to second and third can be obtained as particular cases of general Hermite-Padé problem.

The problem of function asymptotic expansion in prescribed system of functions requires information about class of function [7]. For the function defined by power series this problem can be solved with the aid of mentioned methods by expansion of approximant resulted. This way is concerned with sequence of intermediate expansions with coefficients depending on quantity of power series coefficients used.

As coefficients of intermediate expansions are not equal to coefficients of asymptotic expansion, the difference between partial sums of intermediate expansion functions and original function generally can not be evaluated with the formulas of asymptotic expansion [7].

Also all mentioned approximations require solving a system of linear algebraic equations what reduces their application.

Some recent results developing considered approach are represented in papers [8,9,10]. One result concerned with asymptotic evaluation of defined class of functions using the expansion itself was obtained in [11].

The main approach of present paper is in searching for some function, which Taylor expansion coefficients in fixed point of variable coincide with coefficients of asymptotic expansion of prescribed series in infinite increasing of variable. Using this way we present a sufficient condition of existence and formulas for the coefficients of asymptotic expansion in negative powers of variable for the function defined by power series.

## Notations

$z$ – complex variable, $z = x + i \cdot y$;
$x$ – real variable;
$x_0$ – real value, $x_0 \in (-\infty, \infty)$;
$\bar{x} = \dfrac{1}{x - x_0 + 1}$;

$$\bar{\bar{x}} = \frac{1}{x - x_0};$$

$$F_n(x) = \left(\frac{x}{1-x}\right)^n;$$

$$F'_n(x) = \frac{1}{(1-x)^n};$$

$c_0, ..., c_n, ...$ – real values;

## Definition

*Asymptotically associated function for the function $f(x)$ defined by series*

$$f(x) = \sum_{n=0}^{\infty} c_n \cdot (x - x_0)^n,$$

*is the function $u(x)$ defined by series*

$$u(x) = \sum_{n=1}^{\infty} c_n^* \cdot x^n,$$

*with coefficients*

$$\begin{cases} c_0^* = c_0 \\ c_1^* = c_1 \\ c_n^* = \sum_{s=0}^{n-1} \binom{n-1}{s} \cdot c_{n-s} = \sum_{s=1}^{n-1} \binom{n-1}{s-1} \cdot c_s \end{cases}.$$

**Sufficient condition of existence of asymptotic expansion**
$$f(x) = \sum_{n=0}^{\infty} \frac{q'_n}{(x - x_0 + 1)^n}, \quad x \to \infty$$

Let us consider the next problem. Let $f(z)$ be some function of a complex variable. We know the coefficients of $f(z)$ expansion in positive powers with real center of convergence $(x = x_0)$:

$$f(x) = \sum_{n=0}^{\infty} c_n \cdot (z - x_0)^n. \tag{1}$$

We need the coefficients of $f(x)$ expansion in negative powers

$$f(x) = \sum_{n=0}^{\infty} \frac{q'_n}{(x - x_0 + 1)^n}, \tag{2}$$

if this expansion exists.

This case expansion (2) is asymptotic expansion of $f(x)$ for $x \to \infty$ [7].

Theorem 1. *Let $f(x)$ be a function analytic in point $(x = x_0)$. Then if asymptotically associated function for $f(x)$ is analytic in point $(x = 1)$, then $f(x)$ can be expanded in series (2) and coefficients of this expansion are defined by equations*

$$\begin{cases} q'_0 = u(1) \\ q'_1 = -\dfrac{d}{dx} \cdot u(x) \bigg|_{x=1} \\ q'_n = (-1)^n \cdot \dfrac{1}{n!} \cdot \dfrac{d^n}{dx^n} \cdot u(x) \bigg|_{x=1} \end{cases} \quad (3)$$

Proof. In accordance with Notations,

$$x = \frac{1}{\bar{x}} + x_0 - 1.$$

Formal expansion of $f(x)$ in series (2) is:

$$f(x) = \sum_{n=0}^{\infty} \frac{q'_n}{(x - x_0 + 1)^n},$$

what can be written as

$$v(\bar{x}) = \sum_{n=0}^{\infty} q'_n \cdot \bar{x}^n, \quad (4)$$

where

$$v(\bar{x}) = f\left(\frac{1}{\bar{x}} + x_0 - 1\right). \quad (5)$$

Then, using Taylor formula, we have

$$q'_n = \frac{1}{n!} \cdot \frac{d^n}{dx^n} \cdot v(x) \bigg|_{x=0} \quad (6)$$

Let us introduce the function
$$u(x) = v(1 - x), \quad (7)$$
then
$$v(x) = u(1 - x).$$

The derivatives of *n*-th term of $v(x)$ and $u(x)$ are connected by equation

$$\frac{d^n v(x)}{dx^n} = (-1)^n \cdot \frac{d^n u(1-x)}{d(1-x)^n}, \quad (8)$$

what can be proved by induction.

In accordance with (6),

$$q'_n = \frac{1}{n!} \cdot (-1)^n \cdot \left. \frac{d^n u(1-x)}{d(1-x)^n} \right|_{x=0}, \tag{9}$$

or

$$q'_n = \frac{1}{n!} \cdot (-1)^n \cdot \left. \frac{d^n u(x)}{d(x)^n} \right|_{x=1}. \tag{10}$$

Let $u(x)$ be asymptotically associated function for $f(x)$. In accordance with Theorem 1 hypothesis $u(x)$ is analytic in point $(x=1)$. Then the expansion exists

$$u(x) = \sum_{n=0}^{\infty} a_n \cdot (x-1)^n,$$

where, in accordance with Taylor formula,

$$a_n = \frac{1}{n!} \cdot \left. \frac{d^n u(x)}{dx^n} \right|_{x=1},$$

or, in accordance with (10),

$$a_n = (-1)^n \cdot q'_n.$$

Then

$$u(x) = \sum_{n=0}^{\infty} (-1)^n \cdot q'_n \cdot (x-1)^n. \tag{11}$$

If $u(x)$ expands in series (11) then series $\sum_{n=0}^{\infty} (-1)^n \cdot q'_n \cdot x^n$ converges inside some radius of convergence $R$. This radius of convergence is some function of absolute values of (11) coefficients. E.g. d'Alembert sign of convergence leads to formula

$$R = \lim_{n \to \infty} \sqrt[k]{\frac{|(-1)^n \cdot q'_n|}{|(-1)^{n+k} \cdot q'_{n+k}|}} = \lim_{n \to \infty} \sqrt[k]{\frac{|q'_n|}{|q'_{n+k}|}}. \tag{12}$$

In accordance with (12) each power series with the same absolute values of coefficients to (12) series converges while $(x < R)$.

Then if $u(x)$ is asymptotically associated function for $f(x)$ then series (4) converges while $(\bar{x} < R)$.

Let us prove that $u(x)$ is asymptotically associated function for $f(x)$.

We have

$$u(x) = \sum_{n=0}^{\infty} c_n^* \cdot x^n, \tag{13}$$

where $u(x)$ is defined only by equation (7), $c_n^*$ are values to define.

In accordance with (5) we can write

$$v(1-x) = f\left(\frac{1}{1-x} + x_0 - 1\right),$$

or

$$v(1-x) = f\left(x_0 + \frac{x}{1-x}\right). \tag{14}$$

Under Theorem 1 hypothesis $f(x)$ is analytic in point $(x = x_0)$, so, in accordance with (1), we can write expansion

$$f\left(x_0 + \frac{x}{1-x}\right) = \sum_{n=0}^{\infty} c_n \cdot \left(\frac{x}{1-x}\right)^n. \tag{15}$$

(15) right part converges while

$$x \in \left(-\frac{R}{1-R}, \frac{R}{1+R}\right), \tag{16}$$

$R$ is radius of convergence for series (1).

In accordance with (7) and (13 – 15) we have

$$\sum_{n=0}^{\infty} c_n^* \cdot x^n = \sum_{n=0}^{\infty} c_n \cdot \left(\frac{x}{1-x}\right)^n. \tag{17}$$

Let us consider coefficients of expansion of the function

$$F_n(x) = \left(\frac{x}{1-x}\right)^n \tag{18}$$

in formal power series

$$F_n(x) = \sum_{k=0}^{\infty} c_{n,k} \cdot x^k. \tag{19}$$

In accordance with Notations,

$$F_n(x) = x^n \cdot F_n'(x).$$

$F_n'(x)$ can be expanded in powers of $x$:

$$F_n'(x) = \sum_{k=0}^{\infty} c_{n,k}' \cdot x^k. \tag{20}$$

Let us consider the coefficients $c_{n,k}'$.

The expansion exists

$$F_1'(x) = \frac{1}{1-x} = \sum_{k=0}^{\infty} x^k. \tag{21}$$

Next identity for integers $n > 0$

$$\frac{d^n}{dx^n} \cdot F_1'(x) = \frac{n!}{(1-x)^{n+1}} = n! \cdot F_{n+1}'(x) \tag{22}$$

can be proved by induction.

Then

$$F_n'(x) = \frac{1}{(n-1)!} \cdot \frac{d^{n-1}}{dx^{n-1}} \cdot F_1'(x). \tag{23}$$

Next identity for integers $n > 0$

$$\frac{d^n}{dx^n} \cdot \sum_{k=0}^{\infty} x^k = \sum_{k=n}^{\infty} \prod_{i=0}^{n-1} (k-i) \cdot x^{k-n} \tag{24}$$

can be proved by induction.

In accordance with (21), (23) and (24),

$$F'_n(x) = \frac{1}{(n-1)!} \cdot \sum_{k=n-1}^{\infty} \prod_{i=0}^{n-2}(k-i) \cdot x^{k-n+1} . \tag{25}$$

Let us consider the coefficients of (25) terms. The identity exists

$$\frac{\prod_{i=0}^{n-2}(k-i)}{(n-1)!} \equiv \binom{k}{n-1} . \tag{26}$$

In accordance with (25) and (26),

$$F'_n(x) = \sum_{k=n-1}^{\infty} \binom{k}{n-1} \cdot x^{k-n+1} . \tag{27}$$

Let us introduce summation index $s = k - n + 1$, then $k = s + n - 1$ and (27) can be transformed to

$$F'_n(x) = \sum_{s=0}^{\infty} \binom{s+n-1}{n-1} \cdot x^s . \tag{28}$$

In accordance with (19), (20) and (28),

$$c'_{n,k} = \binom{k+n-1}{n-1}, \tag{29}$$

and

$$F_n(x) = \sum_{k=0}^{\infty} \binom{k+n-1}{n-1} \cdot x^{n+k} . \tag{30}$$

Let us introduce summation index $s = k + n$, then $k = s - n$ and (30) can be transformed to

$$F_n(x) = \sum_{s=n}^{\infty} \binom{s-1}{n-1} \cdot x^s ,$$

then

$$c_{n,k} = \begin{cases} 0 \text{ if } k < n \\ \binom{k-1}{n-1} \text{ if } k \geq n \end{cases},$$

or, in accordance with property of binomial coefficients,

$$c_{n,k} = (-1)^{k-n} \cdot \binom{k-1}{n-1} . \tag{31}$$

Let us equate the coefficients of equal terms of (17) left part expansion in power series and (17) right part and obtain next formulas

$$\begin{cases} c_0^* = c_0 \\ c_1^* = c_1 \\ c_n^* = \sum_{k=1}^{n-1} \binom{n-1}{k-1} \cdot c_k \end{cases} . \tag{33}$$

After replacement in (33) $k$ by $s$ we will obtain equations of Definition for asymptotically associated function.

**On coefficients of asymptotic expansion** $f(x)=\sum_{n=0}^{\infty}\dfrac{q_n}{(x-x_0)^n}$, $x\to\infty$

Let us consider the problem of obtaining the coefficients of expansion

$$f(x)=\sum_{n=0}^{\infty}\frac{q_n}{(x-x_0)^n} \qquad (34)$$

under condition of existence of expansion (2).

If expansion (34) exists then it is asymptotic expansion of $f(x)$ for $x\to\infty$ [7].

Let us express coefficients of (34) in coefficients of expansion (2). In accordance with Notations, functions $\bar{x}$ and $\bar{\bar{x}}$ are connected by the equation

$$\bar{x}=\frac{\bar{\bar{x}}}{1+\bar{\bar{x}}}. \qquad (35)$$

Then formal expansions of $f(x)$ in (2) and (34) can be written as

$$f(x)=\sum_{n=0}^{\infty}q'_n\cdot\bar{x}^n=\sum_{n=0}^{\infty}q_n\cdot\bar{\bar{x}}^n. \qquad (36)$$

In accordance with (35),

$$\sum_{n=0}^{\infty}q'_n\cdot\left(\frac{\bar{\bar{x}}}{1+\bar{\bar{x}}}\right)^n=\sum_{n=0}^{\infty}q_n\cdot\bar{\bar{x}}^n. \qquad (37)$$

The equation (37) can be transformed to (17) type after replacement $\bar{\bar{x}}$ by $(-x)$, $q'_n$ by $\left((-1)^n\cdot q'_n\right)$ and $q_n$ by $\left((-1)^n\cdot q'_n\right)$:

$$\sum_{n=0}^{\infty}(-1)^n\cdot q'_n\cdot\left(\frac{x}{1-x}\right)^n=\sum_{n=0}^{\infty}(-1)^n\cdot q_n\cdot x^n. \qquad (38)$$

Then (33) can be transformed to

$$\begin{cases} q_0=q'_0 \\ -q_1=-q'_1 \\ (-1)^n\cdot q_n=\sum_{k=1}^{n-1}\binom{n-1}{k-1}\cdot(-1)^k\cdot q'_k \end{cases}, \qquad (39)$$

or

$$\begin{cases} q_0=q'_0 \\ q_1=q'_1 \\ q_n=\sum_{k=1}^{n-1}(-1)^{n+k}\cdot\binom{n-1}{k-1}\cdot q'_k \end{cases}. \qquad (40)$$

Formulas (40) uniquely and finitely define coefficients of expansion (34) on prescribed coefficients of expansion (2), so Theorem 1 formulates sufficient conditions of existence for expansion (34).

## Example of method application

Define first two coefficients of expansion (2) for the function $f(x) = \mathrm{arctg}(x)$ using only coefficients of $f(x)$ expansion (1) with center of convergence $(x_0 = 0)$.

In accordance with Theorem 1,

$$\begin{cases} q'_0 = u(1) \\ q'_1 = -\dfrac{du(x)}{dx}\bigg|_{x=1} \end{cases}. \tag{41}$$

Let us calculate these values.

$f(x)$ expands in series (1) as

$$f(x) = \sum_{n=0}^{\infty} \frac{(-1)^n}{2 \cdot n + 1} \cdot x^n,$$

or

$$f(x) = \sum_{n=1}^{\infty} \frac{x^n}{n} \cdot \sin \frac{\pi \cdot n}{2}, \tag{42}$$

so

$$c_n = \frac{1}{n} \cdot \sin \frac{\pi \cdot n}{2}, \tag{43}$$

what can be immediately used for calculation of $c_n^*$ with help of equations (33).

First 30 coefficients of series (42) and (13) with 10 decimal places are represented in Table 1.

Table 1 – coefficients of $f(x)$ and $u(x)$ for $f(x) = \text{arctg}(x)$.

| $n$ | $c_n$ | $c_n^*$ | $n$ | $c_n$ | $c_n^*$ |
|---|---|---|---|---|---|
| 0 | 0 | 0 | 16 | 0 | 0 |
| 1 | 1 | 1 | 17 | 1/17 | 15,05882353 |
| 2 | 0 | 1 | 18 | 0 | 28,44444444 |
| 3 | −1/3 | 0,6666666666 | 19 | −1/19 | 26,94736842 |
| 4 | 0 | 0 | 20 | 0 | 0 |
| 5 | 1/5 | −0,8 | 21 | 1/21 | −48,76190476 |
| 6 | 0 | −1,333333333 | 22 | 0 | −93,09090909 |
| 7 | −1/7 | −1,142857143 | 23 | −1/23 | −89,04347826 |
| 8 | 0 | 0 | 24 | 0 | 0 |
| 9 | 1/9 | 1,777777777 | 25 | 1/25 | 163,84 |
| 10 | 0 | 3,2 | 26 | 0 | 315,0769231 |
| 11 | −1/11 | 2,9090909090 | 27 | −1/27 | 303,4074074 |
| 12 | 0 | 0 | 28 | 0 | 0 |
| 13 | 1/13 | −4,923076923 | 29 | 1/29 | −564,9655172 |
| 14 | 0 | −9,142857143 | 30 | 0 | −1092,262626 |
| 15 | −1/15 | −8,533333333 | 31 | 1/31 | −1057,032258 |

We have discovered that $c_n^*$ values from Table 1 are described by formula:

$$c_n^* = \begin{cases} (-1)^{(n\,\text{div}\,4)} \cdot \dfrac{2^{(n\,\text{div}\,2)}}{n} & \text{if } n \bmod 4 \neq 0 \\ 0 & \text{if } n \bmod 4 = 0 \end{cases}. \qquad (44)$$

In accordance with (12) for $(k = 4)$ and (44), radius of convergence of $u(x)$ is

$$R = \lim_{n \to \infty} \sqrt[4]{\dfrac{2^{(n\,\text{div}\,2)}}{n} \cdot \dfrac{n+4}{2^{((n+4)\,\text{div}\,2)}}} = \lim_{n \to \infty} \sqrt[4]{\dfrac{n+4}{n} \cdot \dfrac{2^{(n\,\text{div}\,2)}}{2^{(n\,\text{div}\,2)+2}}},$$

or

$$R = \dfrac{1}{\sqrt[4]{4}} = \dfrac{1}{\sqrt{2}} \approx 0{,}707106781. \qquad (45)$$

In accordance with (45) $u(x)$ can not be calculated in point $(x = 1)$ by immediate summation of series (13) and we do not know if $u(x)$ is analytic in this point.

To research $u(x)$ let us consider analytic continuations of $u(x)$. Coefficients of analytic continuations can be calculated by general formula:

$$\begin{cases} c_0(x + \Delta x) = \sum_{n=0}^{\infty} c_n(x) \cdot \Delta x^n \\ c_k(x + \Delta x) = \sum_{n=k}^{\infty} c_n(x) \cdot \binom{n}{k} \cdot \Delta x^{n-k} \end{cases}, \qquad (46)$$

where $c_k(x + \Delta x)$ is $k$-th term coefficient of expansion

$$u(a) = \sum_{n=0}^{\infty} c_n(x + \Delta x) \cdot (a - x - \Delta x)^n, \qquad (47)$$

$c_k(x)$ is $k$-th term coefficient of expansion

$$u(a) = \sum_{n=0}^{\infty} c_n(x) \cdot (a - x)^n. \qquad (48)$$

To practically use (46) is necessary to terminate summation after some $n$ fulfilling the inequality

$$c_n(x) \cdot \binom{n}{k} \cdot \Delta x^{n-k} < \alpha, \qquad (49)$$

where $\alpha$ is some prescribed value.

Let us choice a step of variation of (13) analytic continuation center of convergence $\Delta x$.

For fixed values of quantity $m$ of (13) coefficients and acceptable accuracy of calculation analytic continuation coefficients we can obtain for each analytic continuation less coefficients than for previous.

Let us compare results of $c_0(1)$ and $c_1(1)$ calculation for variation of $m$, $\Delta x$ and $\alpha$. Calculation was performed with 19 decimal places. Results are represented in Table 2.

Table 2. Calculated values of $c_0(1)$ and $c_1(1)$ for different values of $m$, $\Delta x$ and $\alpha$.

| $m$ | $\Delta x = 0{,}125$ $\alpha = 0{,}1$ | $\Delta x = 0{,}25$ $\alpha = 0{,}1$ | $\Delta x = 0{,}50$, $\alpha = 0{,}1$ | $\Delta x = 0{,}25$ $\alpha = 0{,}01$ |
|---|---|---|---|---|
| 98 | 1,460166095 – | 1,510539335 0,03453505191 | 1,620321766 – | 1,647310773 – |
| 201 | 1,588720015 – | 1,563427338 0,8928224346 | 1,598489603 1,999998434 | 1,603315576 – |
| 301 | – – | 1,592952061 1,965629445 | 1,552961818 – | 1,563720630 – |
| 401 | 1,569690021 0,9749395282 | 1,553687883 0,05661154018 | 1,556918109 –2,581297548E–7 | 1,560597154 – |
| 501 | 1,569494846 0,9751221205 | 1,570917536 1,003878096 | 1,559558920 0,010434639393 | 1,570410686 0,9868629798 |
| 601 | 1,571060707 0,9896974884 | 1,570813848 1,000744658 | 1,346950822 – | 1,570761086 0,9985219333 |
| 701 | 1,570824591 0,9988651153 | 1,570791841 0,9997836770 | –313,1570000 – | 1,570811918 1,000485872 |
| 801 | 1,570503355 0,9926329933 | 1,570802540 1,000277328 | –807157,0000 – | 1,570771038 0,9990169602 |
| 901 | 1,570411470 0,9910297507 | 1,570770269 0,9989899576 | –3,030051678E10 – | 1,570811829 1,000485872 |
| 1001 | 1,570860320 1,002399327 | 1,570767283 0,9988848701 | –1,983251723E13 – | 1,570811262 1,000468103 |

As most values of $u(x)$ and $\dfrac{du(x)}{dx}$ in point $(x=1{,}25)$ (not represented) calculated by different variants of analytic continuation of series (44) coincide in first decimal places, $u(x)$ seems to be analytic in point $(x=1)$, and in accordance with Theorem 1, values of Table 2 can be used for determination of the coefficients required.

In accordance with (41),
$$q'_0 = c_0(1), \quad q'_1 = -c_1(1). \tag{50}$$

Let us obtain exact values of coefficients using properties of function $f(x) = \operatorname{arctg}(x)$. The identity exists
$$\operatorname{arctg}(x) = \frac{\pi}{2} - \operatorname{arctg}\left(\frac{1}{x}\right). \tag{51}$$

After expansion of (51) right part in series (42) we have
$$\operatorname{arctg}(x) = \frac{\pi}{2} - \sum_{n=1}^{\infty} \frac{1}{n \cdot x^n} \cdot \sin\frac{\pi \cdot n}{2}. \tag{52}$$

Then $q_0 = \dfrac{\pi}{2}$, $q_1 = -1$.

In accordance with (33),
$$q'_0 = \frac{\pi}{2}, \quad q'_1 = -1. \tag{53}$$

In accordance with (50) and (53), exact values of calculated coefficients are
$$c_0(1) = \frac{\pi}{2} \approx 1{,}570796327, \quad c_1(1) = 1.$$

### On coefficients of expansion (3) if (13) converges for $|x| > 1$.

In accordance with the example represented the accuracy of the method depends on several parameters of computing scheme. Not concerning the problem of accuracy increasing let us consider the case when original series (1) allows to obtain formulas for the coefficients of expansion (2) as a converging infinite sum. In accordance with Theorem 1, a sufficient condition of expansion (2) existence for series (1) is analyticity of asymptotically associated function $u(x)$ (13) in point $(x=1)$. Under this condition coefficients of series (2) are defined by equations (3). If radius of convergence of series (13) $R > 1$ coefficients of (2) can be obtained by immediate summation of (13) and its derivatives. To obtain these formulas first let us formulate and prove several propositions concerning binomial coefficients then let us transform equations (3).

### Propositions

Proposition 1. *For integers a and elements of some countable set $r_i$ we have*
$$\sum_{k=0}^{m} \binom{m}{k} \cdot r_{m+a-k} = \sum_{k=0}^{m} \binom{m}{k} \cdot r_{k+a}. \tag{54}$$

Proposition 2. *For integers $m \geq 1$ and $k \geq 0$, $m > k$, we have*

$$\sum_{n=k}^{m}\binom{n-1}{k-1}=\binom{m}{k}. \tag{55}$$

Proposition 3. *For integers $m \geq 1$ and $k \geq 0$ we have*

$$\sum_{z=0}^{m}\binom{k+z}{k}=\binom{k+m+1}{k+1}. \tag{56}$$

Proposition 4. *For integers $m \geq 1, k \geq 1, s \geq 1, r \in [1,k]$, we have*

$$\sum_{n=s}^{m}\binom{n}{k}\cdot\binom{n-1}{s-1}=\sum_{z=0}^{r-1}(-1)^{z}\cdot\binom{m-z}{k-z}\cdot\binom{m}{s+z}+$$
$$+(-1)^{r}\cdot\sum_{n=s}^{m-r}\binom{n}{k-r}\cdot\binom{n-1+r}{s-1+r}. \tag{57}$$

### Proofs of the Propositions

Proposition 1. Let us introduce summation index $z = m - k$, then $k = m - z$. Sum of (54) left part can be written as

$$S_1 = \sum_{m-z=0}^{m}\binom{m}{m-z}\cdot r_{z+a},$$

or

$$S_1 = \sum_{m-z=0}^{m}\binom{m}{z}\cdot r_{z+a}. \tag{58}$$

After reversing of (58) summation direction from $z = m,...,0$ to $z = 0,...,m$, we have

$$S_1 = \sum_{z=0}^{m}\binom{m}{z}\cdot r_{z+a}. \tag{59}$$

After replacement in (59) $z$ with $k$ we will obtain the sum of (54) right part.

Proposition 2. Let us use induction. For $(k = 1)$ right part of (61) is

$$\binom{m}{1} \equiv m,$$

left part

$$\sum_{n=1}^{m}\binom{n-1}{0}=\sum_{n=1}^{m}1=m. \tag{60}$$

For $(k = 1)$ equation (55) is true.

For $(m = 1)$ equation (55) is true.

Let us suppose equation (55) is true for some integer $m$.

Then, if Proposition 2 is true, an equation should be true for $(m + 1)$

$$\sum_{n=k}^{m+1}\binom{n-1}{k-1}=\binom{m+1}{k}. \tag{61}$$

Left part of (61) is
$$\sum_{n=k}^{m+1}\binom{n-1}{k-1} = \sum_{n=k}^{m}\binom{n-1}{k-1}+\binom{m}{k-1}.$$
As equation (55) is true for $m$ by assumption, last sum can be transformed to
$$\sum_{n=k}^{m+1}\binom{n-1}{k-1} = \binom{m}{k}+\binom{m}{k-1}. \tag{62}$$
In accordance with property of binomial coefficients (Pascal's triangle), the identity exists
$$\binom{m}{k}+\binom{m}{k-1} \equiv \binom{m+1}{k}.$$
what coincides with right part of (61).

Proposition 3. Let us use induction. For $(m=1)$ we have
$$\sum_{z=0}^{1}\binom{k+z}{k} = \binom{k}{k}+\binom{k+1}{k},$$
or, as $\binom{k}{k} \equiv 1$ and $\binom{k+1}{k} \equiv k+1$,
$$\sum_{z=0}^{1}\binom{k+z}{k} = k+2. \tag{63}$$
(56) right part can be transformed to
$$\binom{k+m+1}{k+1} = \binom{k+2}{k+1} \equiv k+2,$$
what coincides with (63) right part. So for $(m=1)$ Proposition 3 is true.

Let us suppose Proposition 3 is true for some positive integer $m$.

If Proposition 3 is true for any $m$ then for $(m+1)$ an equation should be true
$$\sum_{z=0}^{m+1}\binom{k+z}{k} = \binom{k+m+2}{k+1}. \tag{64}$$
Left part of (64) can be transformed to
$$\sum_{z=0}^{m+1}\binom{k+z}{k} = \sum_{z=0}^{m}\binom{k+z}{k}+\binom{k+m+1}{k}.$$
As (56) is true by supposition, the equation is true
$$\sum_{z=0}^{m+1}\binom{k+z}{k} = \binom{k+m+1}{k+1}+\binom{k+m+1}{k}.$$
The identity exist [12] (Pascal's triangle):
$$\binom{a}{b}+\binom{a}{b+1} \equiv \binom{a+1}{b+1},$$
then
$$\binom{k+m+1}{k+1}+\binom{k+m+1}{k} = \binom{k+m+2}{k+1}. \tag{65}$$

what coincides with right part of (64).

Proposition 4. Let us use induction. If (57) is true for $(r=1)$ then the equation should be true

$$\sum_{n=s}^{m}\binom{n}{k}\cdot\binom{n-1}{s-1}=\binom{m}{k}\cdot\binom{m}{s}-\sum_{n=s}^{m-1}\binom{n}{k-1}\cdot\binom{n}{s}. \quad (66)$$

Let us introduce summation index $p=n-s+1$, then $n=p+s-1$, for $n=s$ $p=1$, for $n=m$ $p=m-s+1$, then

$$\sum_{n=s}^{m}\binom{n}{k}\cdot\binom{n-1}{s-1}=\sum_{p=1}^{m-s+1}\binom{p+s-1}{k}\cdot\binom{p+s-2}{s-1}.$$

To transform left part of (66) let us use Abel's lemma [7]:

$$\sum_{n=s}^{m}\binom{n}{k}\cdot\binom{n-1}{s-1}=\sum_{p=1}^{m-s}A_p\cdot(b_p-b_{p+1})+A_{m-s+1}\cdot b_{m-s+1}, \quad (67)$$

where

$$A_p=\sum_{l=1}^{p}\binom{l+s-2}{s-1},\quad A_{m-s+1}=\sum_{l=1}^{m-s+1}\binom{l+s-2}{s-1},$$

$$b_p=\binom{p+s-1}{k},\ b_{p+1}=\binom{p+s}{k},\ b_p-b_{p+1}=-\binom{p+s-1}{k-1},$$

$$b_{m-s+1}=\binom{m}{k}.$$

We have

$$A_p=\binom{p+s-1}{s},$$

$$A_{m-s+1}=\binom{m}{s}.$$

Then

$$\sum_{n=s}^{m}\binom{n}{k}\cdot\binom{n-1}{s-1}=\binom{m}{k}\cdot\binom{m}{s}-\sum_{p=1}^{m-s}\binom{p+s-1}{s}\cdot\binom{p+s-1}{k-1}. \quad (68)$$

After returning to summation index $n$ in the second additive of (68) right part,

$$\sum_{n=s}^{m}\binom{n}{k}\cdot\binom{n-1}{s-1}=\binom{m}{k}\cdot\binom{m}{s}-\sum_{n=s}^{m-1}\binom{n}{s}\cdot\binom{n}{k-1},$$

what coincides with (66).

Then for $(r=1)$ equation (57) is true.

Let (57) be true for integer $a$, $1<a<k$.

If Proposition 4 is true then the equation should be true

$$\sum_{n=s}^{m}\binom{n}{k}\cdot\binom{n-1}{s-1}=\sum_{z=0}^{r}(-1)^z\cdot\binom{m-z}{k-z}\cdot\binom{m}{s+z}+$$

$$(-1)^{r+1}\cdot\sum_{n=s}^{m-r-1}\binom{n}{k-r-1}\cdot\binom{n+r}{s+r}. \quad (69)$$

Let us introduce summation index $p = n - s + 1$ then $n = p + s - 1$, for $n = s$ $p = 1$, for $n = m - r$ $p = m - r - s + 1$, then second additive of (9) right part is

$$S = (-1)^r \cdot \sum_{p=1}^{m-r-s+1} \binom{p+s-1}{k-r} \cdot \binom{p+s+r-2}{s-1+r}. \qquad (70)$$

To transform (70) let us use Abel's lemma:

$$S = (-1)^r \cdot \left\{ \sum_{p=1}^{m-r-s} A_p \cdot (b_p - b_{p+1}) + A_{m-r-s+1} \cdot b_{m-r-s+1} \right\}, \qquad (71)$$

$$A_p = \sum_{l=1}^{p} \binom{l+s+r-2}{s-1+r}, \quad A_{m-r-s+1} = \sum_{l=1}^{m-r-s+1} \binom{l+r+s-2}{s-1+r},$$

$$b_p = \binom{p+s-1}{k-r}, b_{p+1} = \binom{p+s}{k-r}, \; b_p - b_{p+1} = -\binom{p+s-1}{k-r-1},$$

$$b_{m-r-s+1} = \binom{m-r}{k-r}.$$

Let us introduce summation index $t = l - 1$, then $l = t + 1$, for $l = 1$ $t = 0$, for $l = p$ $t = p - 1$, then

$$A_p = \sum_{t=0}^{p-1} \binom{t+s+r-1}{s+r-1}. \qquad (72)$$

In accordance with (56),

$$A_p = \binom{s+r-1+p-1+1}{s+r-1+1} = \binom{r+p+s-1}{s+r}, \qquad (73)$$

in particulars,

$$A_{m-r-s+1} = \binom{r+m-r-s+1+s-1}{s+r} = \binom{m}{s+r}. \qquad (74)$$

In accordance with (71 – 74),

$$S = (-1)^r \cdot \left\{ \binom{m}{s+r} \cdot \binom{m-r}{k-r} - \sum_{p=1}^{m-r-s} \binom{r+p+s-1}{s+r} \cdot \binom{p+s-1}{k-r-1} \right\}. \qquad (75)$$

After returning to summation index $n$ in the second additive of (75) right part,

$$S = (-1)^r \cdot \left\{ \binom{m}{s+r} \cdot \binom{m-r}{k-r} - \sum_{n=s}^{m-r-1} \binom{r+n}{r+s} \cdot \binom{n}{k-r-1} \right\}. \qquad (76)$$

As the second additive of (57) right part is equal (76), (57) can be transformed to

$$\sum_{n=s}^{m} \binom{n}{k} \cdot \binom{n-1}{s-1} = \sum_{z=0}^{r-1} (-1)^z \cdot \binom{m-z}{k-z} \cdot \binom{m}{s+z} +$$

$$+ (-1)^r \cdot \binom{m-r}{k-r} \cdot \binom{m}{s+r} + (-1)^{r+1} \cdot \sum_{n=s}^{m-r-1} \binom{n}{k-r-1} \cdot \binom{n+r}{s+r}. \qquad (77)$$

The second additive of (77) right part can be put under the summation operator of the first additive for $(z = r)$ then (77) coincides with (69).

## Transformation of equations (3)

Derivatives of a converging power series are described with identity

$$\frac{1}{k!} \cdot \frac{d^k}{dx^k} \cdot \sum_{n=0}^{\infty} a_n \cdot x^n = \sum_{n=k}^{\infty} \binom{n}{k} \cdot a_n \cdot x^{n-k}, \qquad (78)$$

$n, k$ – integer values, $x$ is inside radius of convergence of (78).

In accordance with (78) equations (3) can be written as

$$\begin{cases} q'_0 = \sum_{n=0}^{\infty} c_n^* \\ q'_k = (-1)^k \cdot \sum_{n=k}^{\infty} \binom{n}{k} \cdot c_n^* \end{cases}, \qquad (79)$$

Let us obtain equations connecting $q'_k$ and $c_n$ where each coefficient of expansion (1) is used only once. Obviously if (33) is directly substituted to (79) we will have instead of required formulas an infinite sum of terms containing a product of finite and infinite values. To eliminate the uncertainty let us consider the converging series as a limit of a finite sum for infinite increasing of processing terms quantity.

As all considered series in (79) are converging, we have

$$\begin{cases} q'_0 = \lim_{m \to \infty} \sum_{n=0}^{m} c_n^* \\ q'_k = (-1)^k \cdot \lim_{m \to \infty} \sum_{n=k}^{m} \binom{n}{k} \cdot c_n^* \end{cases}. \qquad (80)$$

After substituting (33) to (80) we will obtain

$$\begin{cases} q'_0 = c_0 + \lim_{m \to \infty} \sum_{n=1}^{m} \sum_{s=0}^{n-1} \binom{n-1}{s} \cdot c_{n-s} \\ q'_k = (-1)^k \cdot \lim_{m \to \infty} \sum_{n=k}^{m} \binom{n}{k} \cdot \sum_{s=0}^{n-1} \binom{n-1}{s} \cdot c_{n-s} \end{cases}. \qquad (81)$$

Let us consider finite sums of equations (81).

### Transformation of sum $c_0 + \sum_{n=1}^{m} \sum_{s=0}^{n-1} \binom{n-1}{s} \cdot c_s$

In accordance with (54), first equation of (81) can be written as

$$q'_0 = c_0 + \lim_{m \to \infty} \sum_{n=1}^{m} \sum_{s=0}^{n-1} \binom{n-1}{s} \cdot c_{s+1}.$$

Let us define the sum $q'_{0,m}$ as

$$q'_{0,m} = c_0 + \sum_{n=1}^{m} \sum_{s=0}^{n-1} \binom{n-1}{s} \cdot c_{s+1}. \qquad (82)$$

Let us introduce summation index $z = s+1$ then $s = z-1$, for $s=0$ $z=1$, for $s = n-1$ $z = n$, then

$$q'_{0,m} = c_0 + \sum_{n=1}^{m}\sum_{z=1}^{n}\binom{n-1}{z-1} \cdot c_z.$$

After replacement $z$ with $s$ we have

$$q'_{0,m} = c_0 + \sum_{n=1}^{m}\sum_{s=1}^{n}\binom{n-1}{s-1} \cdot c_s. \tag{83}$$

Let us change the order of (83) summation. As in original sum the coefficient of $c_s$ is $\binom{n-1}{s-1}$ and formula (83) summarizes $c_s$ from $s=1$ to $s=m$, we have

$$q'_{0,m} = c_0 + \sum_{s=1}^{m} c_s \cdot \sum_{n=s}^{m}\binom{n-1}{s-1}. \tag{84}$$

In accordance with (7) let us modify (36) to

$$q'_{0,m} = c_0 + \sum_{s=1}^{m} c_s \cdot \binom{m}{s},$$

or, as $\binom{m}{0} \equiv 1$,

$$q'_{0,m} = \sum_{s=0}^{m} c_s \cdot \binom{m}{s}. \tag{85}$$

**Transformation of sum** $S_k = \sum_{n=k}^{m}\binom{n}{k} \cdot \sum_{s=0}^{n-1}\binom{n-1}{s} \cdot c_{n-s}$

The original equation is

$$q'_{k,m} = (-1)^k \cdot \sum_{n=k}^{m}\binom{n}{k} \cdot \sum_{s=0}^{n-1}\binom{n-1}{s} \cdot c_{n-s}. \tag{86}$$

In accordance with (54) (86) can be modified to

$$q'_{k,m} = (-1)^k \cdot \sum_{n=k}^{m}\binom{n}{k} \cdot \sum_{s=1}^{n}\binom{n-1}{s-1} \cdot c_{s+1}. \tag{87}$$

Let us change the order of (87) summation. Similarly to obtaining of (84) we can write

$$q'_{k,m} = (-1)^k \cdot \sum_{n=k}^{m}\binom{n}{k} \cdot \sum_{s=1}^{n}\binom{n-1}{s-1} \cdot c_s, \tag{88}$$

or

$$q'_{k,m} = (-1)^k \cdot \left\{\binom{k}{k} \cdot \sum_{s=1}^{k}\binom{k-1}{s-1} \cdot c_s + \binom{k+1}{k} \cdot \sum_{s=1}^{k+1}\binom{k}{s-1} \cdot c_s + ... + \right.$$

$$+\binom{k+2}{k} \cdot \sum_{s=1}^{k+2}\binom{k+1}{s-1} \cdot c_s + \ldots + \binom{m}{k} \cdot \sum_{s=1}^{m}\binom{m-1}{s-1} \cdot c_s . \tag{89}$$

or

$$q'_{k,m} = (-1)^k \cdot \sum_{n=1}^{m} \gamma_n \cdot c_n , \tag{90}$$

where $\gamma_n$ can be written in two different ways for $n \leq k$ and for $n > k$:

$$\begin{cases} \gamma_n = \sum_{z=k}^{m}\binom{z}{k} \cdot \binom{z-1}{n-1} & \text{for } n \leq k \\ \gamma_n = \sum_{z=n}^{m}\binom{z}{k} \cdot \binom{z-1}{n-1} & \text{for } n > k \end{cases} . \tag{91}$$

As binomial coefficient of power below zero is equal zero, second equation of (91) is automatically obtained from first equation for $n > k$.

So unitary equation is

$$q'_{k,m} = (-1)^k \cdot \sum_{s=1}^{m} c_s \cdot \sum_{n=k}^{m}\binom{n}{k} \cdot \binom{n-1}{s-1}. \tag{92}$$

Similarly, for $n \leq k$ first equation of (91) is automatically obtained from second equation so another form of unitary equation is

$$q'_{k,m} = (-1)^k \cdot \sum_{s=1}^{m} c_s \cdot \sum_{n=s}^{m}\binom{n}{k} \cdot \binom{n-1}{s-1}. \tag{93}$$

So both equations (92) and (93) are equivalent.

Let us consider $q'_{k,m}$ for different $k$ to obtain unitary formula for $q'_{k,m}$ not containing $m$ summation operations.

**Transformation of sum** $S_k = \sum_{n=k}^{m}\binom{n}{k} \cdot \sum_{s=0}^{n-1}\binom{n-1}{s} \cdot c_{n-s}$ **for** $(k=1)$

Let us use equation (93):

$$q'_{1,m} = -\sum_{s=1}^{m} c_s \cdot \sum_{n=s}^{m}\binom{n}{1} \cdot \binom{n-1}{s-1}. \tag{94}$$

Let us consider sum of (94) right part

$$S = \sum_{n=s}^{m}\binom{n}{1} \cdot \binom{n-1}{s-1} = \sum_{n=s}^{m}\binom{n-1}{s-1} \cdot n . \tag{95}$$

To transform (95) let us use Abel's lemma

$$\sum_{n=1}^{m} a_n \cdot b_n = \sum_{n=1}^{m-1} A_n \cdot (b_n - b_{n+1}) + A_m \cdot b_m , \tag{96}$$

where

$$A_n = \sum_{k=1}^{n} a_k .$$

Let us introduce summation index $z = n - s + 1$, then $n = z + s - 1$, for $n = s$ $z = 1$, for $n = m$ $z = m - s + 1$, so that

$$S = \sum_{z=1}^{m-s+1} (z + s - 1) \cdot \binom{z + s - 2}{s - 1}, \tag{97}$$

or

$$S = \sum_{z=1}^{m-s} A_z \cdot (b_z - b_{z+1}) + A_{m-s+1} \cdot b_{m-s+1}, \tag{98}$$

where

$$A_z = \sum_{k=1}^{z} \binom{k + s - 2}{s - 1}, \qquad A_{m-s+1} = \sum_{k=1}^{m-s+1} \binom{k + s - 2}{s - 1},$$

$b_z = z + s - 1$, $b_{z+1} = z + s$, $b_z - b_{z+1} = -1$, $b_{m-s+1} = m$.

Let us introduce summation index $n = k + s - 1$, then $k = n - s + 1$, for $k = 1$ $n = s$, for $k = z$ $n = z + s - 1$.

Then

$$A_z = \sum_{n=s}^{z+s-1} \binom{n - 1}{s - 1}.$$

In accordance with (55),

$$A_z = \binom{z + s - 1}{z}, \tag{99}$$

in particulars,

$$A_{m-s+1} = \binom{m}{s}. \tag{100}$$

In accordance with (99) and (100),

$$S = m \cdot \binom{m}{s} - \sum_{z=1}^{m-s} \binom{z + s - 1}{s}. \tag{101}$$

Let us introduce summation index $r = z - 1$, then $z = r + 1$, for $z = 1$ $r = 0$, for $z = m - s$ $r = m - s - 1$, so (101) right part can be transformed to

$$\sum_{z=1}^{m-s} \binom{z + s - 1}{s} = \sum_{z=0}^{m-s-1} \binom{z + s}{s}.$$

In accordance with (56),

$$\sum_{z=0}^{m-s-1} \binom{z + s}{s} = \binom{m}{s + 1}. \tag{102}$$

In accordance with (94), (95), (101) and (102),

$$q'_{1,m} = -\sum_{s=1}^{m} c_s \cdot \left\{ m \cdot \binom{m}{s} - \binom{m}{s + 1} \right\}. \tag{103}$$

**Transformation of sum** $S_k = \sum_{n=k}^{m} \binom{n}{k} \cdot \sum_{s=0}^{n-1} \binom{n - 1}{s} \cdot c_{n-s}$ **for** $(k = 2)$

Let us use equation (93):
$$q'_{2,m} = \sum_{s=1}^{m} c_s \cdot \sum_{n=s}^{m} \binom{n}{2} \cdot \binom{n-1}{s-1}. \qquad (104)$$

Let us consider sum of (56) right part:
$$S = \sum_{n=s}^{m} \binom{n}{2} \cdot \binom{n-1}{s-1}. \qquad (105)$$

To transform (105) let us use Abel's lemma.

Let us introduce summation index $z = n - s + 1$, then $n = z - s + 1$, for $n = s$ $z = 1$, for $n = m$ $z = m - s + 1$, then
$$S = \sum_{z=1}^{m-s+1} \binom{z+s-1}{2} \cdot \binom{z+s-2}{s-1}. \qquad (106)$$

In accordance with (96) (106) can be written as
$$S = \sum_{z=1}^{m-s} A_z \cdot (b_z - b_{z+1}) + A_{m-s+1} \cdot b_{m-s+1}, \qquad (107)$$

where
$$A_z = \sum_{k=1}^{z} \binom{k+s-2}{s-1}, \quad A_{m-s+1} = \sum_{k=1}^{m-s+1} \binom{k+s-2}{s-1},$$
$$b_z = \binom{z+s-1}{2}, \ b_{z+1} = \binom{z+s}{2}, \ b_z - b_{z+1} = -\binom{z+s-1}{1} \equiv -(z+s-1),$$
$$b_{m-s+1} = \binom{m}{2}.$$

In accordance with (99) and (100),
$$A_z = \binom{z+s-1}{s},$$
$$A_{m-s+1} = \binom{m}{s}.$$

Then (107) can be transformed to
$$S = \binom{m}{2} \cdot \binom{m}{s} - \sum_{z=1}^{m-s} (z+s-1) \cdot \binom{z+s-1}{s}. \qquad (108)$$

Let us consider the second additive of (108) without its sign:
$$S_2 = \sum_{z=1}^{m-s} (z+s-1) \cdot \binom{z+s-1}{s}. \qquad (109)$$

In accordance with (96) (109) can be written as
$$S_2 = \sum_{z=1}^{m-s-1} A_z \cdot (b_z - b_{z+1}) + A_{m-s} \cdot b_{m-s}, \qquad (110)$$

where

$$A_z = \sum_{k=1}^{z} \binom{k+s-1}{s}, \quad A_{m-s} = \sum_{k=1}^{m-s} \binom{k+s-1}{s},$$

$b_z = z+s-1$, $b_{z+1} = z+s$, $b_z - b_{z+1} = -1$, $b_{m-s} = m-1$.

Let us introduce summation index $r = k-1$ then $k = r+1$, for $k=1$ $r=0$, for $k = z$ $r = z-1$, so that

$$A_z = \sum_{r=0}^{z-1} \binom{r+s}{s}.$$

In accordance with (102),

$$A_z = \binom{s+z}{s+1}, \tag{111}$$

then

$$A_{m-s} = \binom{m}{s+1}. \tag{112}$$

In accordance with (111) and (112), (110) can be transformed to

$$S_2 = (m-1) \cdot \binom{m}{s+1} - \sum_{z=1}^{m-s-1} \binom{s+z}{s+1}. \tag{113}$$

Let us introduce summation index $r = z-1$ then $z = r+1$, for $z=1$ $r=0$, for $z = m-s-1$ $r = m-s-2$, then

$$\sum_{z=1}^{m-s-1} \binom{s+z}{s+2} = \sum_{r=0}^{m+s-2} \binom{s+1+r}{s+1}. \tag{114}$$

In accordance with (104),

$$\sum_{r=0}^{m-s-2} \binom{s+1+r}{s+1} = \binom{m}{s+2}. \tag{115}$$

Then

$$S_2 = (m-1) \cdot \binom{m}{s+1} - \binom{m}{s+2}, \tag{116}$$

$$S = \binom{m}{2} \cdot \binom{m}{s} - \binom{m-1}{1} \cdot \binom{m}{s+1} + \binom{m}{s+2} \cdot \binom{m-2}{0},$$

$$q'_{2,m} = \sum_{s=1}^{m} c_s \cdot \left\{ \binom{m}{2} \cdot \binom{m}{s} - \binom{m-1}{1} \cdot \binom{m}{s+1} + \binom{m-2}{0} \cdot \binom{m}{s+2} \right\}. \tag{117}$$

(55) can be transformed to

$$q'_{1,m} = -\sum_{s=1}^{m} c_s \cdot \left\{ \binom{m}{1} \cdot \binom{m}{s} - \binom{m-1}{0} \right\}. \tag{118}$$

Then (117) and (118) can be written unitary:

$$q'_{k,m} = (-1)^k \cdot \sum_{s=1}^{m} c_s \cdot \sum_{n=0}^{k} (-1)^n \cdot \binom{m-n}{k-n} \binom{m}{s+n}, \tag{119}$$

what is true for $k = 1, ..., 2$.

Let us prove that (119) is true for any integer positive $k$.

If $r$ in (57) is equal $k$ then
$$\sum_{n=s}^{m}\binom{n-1}{s-1}\binom{n}{k} = \sum_{z=0}^{k-1}(-1)^z \cdot \binom{m-z}{k-z}\binom{m}{s+z} +$$
$$+(-1)^k \cdot \sum_{n=s}^{m-k}\binom{n}{0}\binom{n-1+k}{s-1+k}. \tag{120}$$

Let us consider the second additive of (120) right part.

As $\binom{n}{0} \equiv 1$, the additive can be transformed to

$$S = (-1)^k \cdot \sum_{n=s}^{m-k}\binom{n-1+k}{s-1+k}. \tag{121}$$

Let us introduce summation index $t = n-s$ then $n = t+s$, for $n=s$ $t=0$, for $n = m-k$ $t = m-k-s$, then

$$S = (-1)^k \cdot \sum_{t=0}^{m-k-s}\binom{t+s+k-1}{s+k-1}. \tag{122}$$

In accordance with (56), (122) can be transformed to

$$S = (-1)^k \cdot \binom{s+k-1+m-k-s+1}{s+k-1+1} = (-1)^k \cdot \binom{m}{s+k}. \tag{123}$$

In accordance with (123), (120) can be transformed to

$$\sum_{n=s}^{m}\binom{n}{k}\binom{n-1}{s-1} = \sum_{z=0}^{k}(-1)^z \cdot \binom{m-z}{k-z}\binom{m}{s+z}. \tag{124}$$

In accordance with (93),

$$q'_{k,m} = (-1)^k \cdot \sum_{s=1}^{m} c_s \cdot \sum_{z=0}^{k}(-1)^z \cdot \binom{m-z}{k-z}\binom{m}{s+z}. \tag{125}$$

After replacement in (125) $z$ by $n$ (125) coincides with (119).
In accordance with (81), (85) and (125), the final equations are

$$\begin{cases} q'_0 = \lim_{m \to \infty} \sum_{s=0}^{m} c_s \cdot \binom{m}{s} \\ q'_k = (-1)^k \cdot \lim_{m \to \infty} \sum_{s=1}^{m} c_s \cdot \sum_{n=0}^{k}(-1)^n \cdot \binom{m-n}{k-n}\binom{m}{s+1} \end{cases}. \tag{126}$$

## Conclusion

In present paper we have formulated and proved Theorem 1 providing sufficient condition of existence of asymptotic expansion

$$f(x) = \sum_{n=0}^{\infty} \frac{q'_n}{(x - x_0 + 1)^n}, \quad x \to \infty$$

of function defined by series

$$f(x) = \sum_{n=0}^{\infty} c_n \cdot (x - x_0)^n.$$

The condition is analyticity of the function $u(x)$, defined uniquely by coefficients of original series with equations (13) and (33), in point $(x = 1)$.

This condition is also sufficient for existence of asymptotic expansion of $f(x)$

$$f(x) = \sum_{n=0}^{\infty} \frac{q_n}{(x - x_0)^n}, \quad x \to \infty.$$

Theorem 1 can be used for immediate calculation of coefficients of expansions (2) and (34) without solving a system of linear algebraic equations.

During testing of our method with the aid of given function $f(x) = \arctg(x)$ we have obtained results depending on parameters of computational scheme.

The best value of zero term coefficient of $f(x)$ asymptotic expansion contains 6 true decimal places, of first term – 4 places, both coefficients were obtained under the same parameters of computational scheme: quantity of given coefficients of asymptotically associated function $m = 701$, step of center of convergence variation for analytic continuations $\Delta x = 0{,}25$, summation is terminated under current coefficient is less than $\alpha = 0{,}1$. We notify that after some hundred of coefficients processing almost all variants of computational scheme give not less than two true decimal places. Exceptions are variants with step $\Delta x = 0{,}5$ where computational scheme becomes unstable in increasing of quantity of coefficients processed. We suppose that accuracy lost in further increasing of $m$ is due to fixed accuracy of values processing with 19 decimal places so a problem of multi-parameter optimization takes place for all used parameters.

The problem of representation of coefficients of asymptotic expansion

$$f(x) = \sum_{n=0}^{\infty} \frac{q'_n}{(x - x_0 + 1)^n}, \quad x \to \infty$$

through the coefficients of original Taylor expansion (1) was solved by means of combinatorics. Using Propositions 1 – 4, induction and Abel's lemma we have obtained resulted equations (126). In accordance with (33) these formulas can be used if radius of convergence $R$ of series (13) $R > 1$.